\documentclass[amstex,12pt, amssymb]{article}

\usepackage{mathtext}
\usepackage[cp1251]{inputenc}
\usepackage[T2A]{fontenc}
\usepackage[dvips]{graphicx}
\usepackage{amsmath}
\usepackage{amssymb}
\usepackage{amsxtra}
\usepackage{latexsym}
\usepackage{ifthen}

\newcounter{lemma}[section]

\newcounter{corollary}[section]

\newcounter{remark}[section]

\newcounter{theorem}[section]

\newcounter{proposition}[section]

\numberwithin{equation}{section}

\def\XXint#1#2#3{{\setbox0=\hbox{$#1{#2#3}{\int}$}
     \vcenter{\hbox{$#2#3$}}\kern-.5\wd0}}

\def\cc{\setcounter{equation}{0}
\setcounter{figure}{0}\setcounter{table}{0}}

\begin{document}

\markboth{\centerline{Vladimir Ryazanov and Artyem Yefimushkin}}
{\centerline{On Riemann-Hilbert problem for Beltrami equations}}

\author{{Vladimir Ryazanov and Artyem Yefimushkin}}

\title{{\bf ON THE RIEMANN-HILBERT PROBLEM\\
FOR THE BELTRAMI EQUATIONS}}

\maketitle

\large \begin{abstract} It is developed the theory of the Dirichlet
problem for harmonic functions. On this basis, for the nondegenerate
Beltrami equations in the quasidisks and, in particular, in smooth
Jordan domains, it is proved the existence of regular solutions of
the Riemann-Hilbert problem with coefficients of bounded variation
and boundary data that are measurable with respect to the absolute
harmonic measure (logarithmic capacity). Moreover, it was shown that
the dimension of the spaces of the given solutions is infinite.
\end{abstract}


\large \cc
\section{Introduction}

Let $D$ be a domain in the complex plane $\mathbb C$ and let $\mu:
D\to\mathbb C$ be a measurable function with $|\mu(z)|<1$ a.e. The
equation of the form
\begin{equation}\label{1}
f_{\bar{z}}=\mu(z)\cdot f_z\
\end{equation}
where $f_{\bar z}={\bar\partial}f=(f_x+if_y)/2$, $f_{z}=\partial
f=(f_x-if_y)/2$, $z=x+iy$, $f_x$ and $f_y$ are partial derivatives
of the function $f$ in $x$ and $y$, respectively, is said to be a
{\bf Beltrami equation}. The Beltrami equation~\eqref{1} is said to
be {\bf nondegenerate} if $||\mu||_{\infty}<1$.

Note that there were recently established a great number of new
theorems on the existence and on the boundary behavior of
homeomorphic solutions and, on this basis, on the Dirichlet problem
for the Beltrami equations with essentially unbounded {\bf
distortion quotients} $K_{\mu}(z)=(1+|\mu(z)|)/(1-|\mu(z)|)$, see,
e.g., the monographs \cite{AIM}--\cite{MRSY} and the papers
\cite{KPR1}--\cite{RSSY}, and many references therein. However,
under the study of the Riemann-Hilbert problem for ~\eqref{1} we
restrict ourselves here with the nondegenerate case because this
investigation leads to a very delicate Lusin's problem on
interconnections of the boundary data of conjugate harmonic
functions and with the dificult problem on the distortion of
boundary measures under more general mappings.

Boundary value problems for analytic functions are due to the
well-known Riemann dissertation (1851), and also to works of Hilbert
(1904, 1912, 1924), and Poincar\'e (1910), see the monograph
\cite{Vek} for details and also for the case of generalized analytic
functions.

The first concrete problem of such a type has been proposed by
Hilbert (1904) and called at present by the Hilbert problem or the
Riemann-Hilbert problem. That consists in finding analytic functions
$f$ in a domain bounded by a rectifiable Jordan curve $C$ with the
linear boundary condition
\begin{equation}\label{2}
\rm {Re}\ \overline{\lambda(\zeta)}\cdot f(\zeta)\ =\ \varphi(\zeta)
\quad\quad\quad\ \ \ \forall \ \zeta\in C
\end{equation}
where it was assumed by him that the functions $\lambda$ and
$\varphi$ are continuously differentiable with respect to the
natural parameter $s$ on $C$ and, moreover, $|\lambda|\ne 0$
everywhere on $C$. Hence without loss of generality one may assume
that $|\lambda|\equiv 1$ on $C$.

The first way for solving this problem based on the theory of
singular integral equations was given by Hilbert, see \cite{H1}.
This attempt was not quite successful because of the theory of
singular integral equations has been not yet enough developed at
that time. However, just that way became the main approach in this
research direction, see e.g. \cite{Vek}, \cite{G} and \cite{Mus}. In
particular, the existence of solutions to this problem was  in that
way proved for H\"older continuous $\lambda$ and $\varphi$, see e.g.
\cite{G}.

Another way for solving this problem based on a reduction to the
corresponding two Dirichlet problems was also proposed by Hilbert,
see e.g. \cite{H2}. A very general solution of the Riemann-Hilbert
problem by this way was  recently given in \cite{Ry} for the
arbitrary Jordan domains with functions $\varphi$ and $\lambda$ that
are only measurable with respect to the harmonic measure.

We follow the second scheme of Hilbert under the study of the
generalized Riemann-Hilbert problem for the Beltrami equations.
However, as it follows from the known Ahlfors--Beurling--Bishop
examples, see e.g. \cite{AB$_2$}, the harmonic measure zero is not
invariant under quasiconformal mappings. Hence here we apply the
so--called absolute harmonic measure (logarithmic capacity).

Recall that homeomorphic solutions with distributional derivatives
of the nondegenerate Beltrami equations~\eqref{1} are called {\bf
quasiconformal mappings}, see e.g. \cite{Alf} and \cite{LV}. The
images of the unit disk $\mathbb D = \{ z\in\mathbb{C}: |z|<1\}$
under the quasiconformal mappings $\mathbb C$ onto itself are called
{\bf quasidisks} and their boundaries are called {\bf quasicircles}
or {\bf quasiconformal curves}. Recall that a {\bf Jordan curve} is
a continuous one-to-one image of the unit circle in $\mathbb C$. It
is known that every smooth (or Lipschitz) Jordan curve is a
quasiconformal curve and, at the same time, quasiconformal curves
can be nonrectifiable as it follows from the known examples, see
e.g. the point II.8.10 in \cite{LV}.

Note that a Jordan curve generally speaking has no tangents. Hence
we need a replacement for the notion of a nontangential limit
usually applied. In this connection, recall the Bagemihl theorem in
\cite{B}, see also Theorem III.1.8 in \cite{No}, stated that, for
any function $\Omega:\mathbb D\to\overline{\mathbb C}$, except at
most countable set of points $\zeta\in\partial\mathbb D$, for all
pairs of arcs $\gamma_1$ and $\gamma_2$ in $\Bbb D$ terminating at
$\zeta\in\partial\mathbb D$,
\begin{equation}\label{3}
C(\Omega,\gamma_1)\ \cap\ C(\Omega,\gamma_2)\ \neq\ \varnothing\ ,
\end{equation}
where $C(\Omega,\gamma)$ denotes the {\bf cluster set of $\Omega$ at
$\zeta$ along $\gamma$}, i.e.,
$$
C(\Omega,\gamma)\ =\ \{ w\in\overline{\mathbb C}\ :\ \Omega(z_n)\to
w,\ z_n\to\zeta ,\ z_n\in\gamma\}\ .
$$
Immediately by the theorems of Riemann and Caratheodory, this result
is extended to an arbitrary Jordan domain $D$ in $\mathbb C$. Given
a function $\Omega: D\to\overline{\mathbb C}$ and $\zeta\in\partial
D$, denote by $P(\Omega , \zeta)$ the intersection of all cluster
sets $C(\Omega,\gamma)$ for arcs $\gamma$ in $D$ terminating at
$\zeta$. Later on, we call the points of the set $P(\Omega , \zeta)$
{\bf principal asymptotic values} of $\Omega$ at $\zeta$. Note that
if $\Omega$ has a limit along at least one arc in $D$ terminating at
a point $\zeta\in\partial D$ with the property (\ref{3}), then the
principal asymptotic value is unique.

Recall also that a mapping $f:{D}\to{\Bbb C}$ is called {\bf
discrete} if the pre-image $f^{-1}(y)$ consists of isolated points
for every $y\in \Bbb C$, and {\bf open} if $f$ maps every open set
$U\subseteq D$ onto an open set in $\Bbb C$.

The regular solution of a Beltrami equation~\eqref{1} is a
continuous, discrete and open mapping $f:{D}\to{\Bbb C}$ with
distributional derivatives satisfying~\eqref{1} a.e. Note that, in
the case of nondegenerate Beltrami equations~\eqref{1}, a regular
solution $f$ belongs to class $W_{\rm{loc}}^{1,p}$ for some $p>2$
and, moreover, its Jacobian $J_f(z)\not\equiv 0$ for almost all
$z\in D$, and it is called a {\bf quasiconformal function}, see e.g.
Chapter VI in \cite{LV}.

A {\bf regular solution of the Riemann-Hilbert problem}~\eqref{2}
for the Beltrami equation~\eqref{1} is a regular solution of
\eqref{1} satisfying the boundary condition \eqref{2} in the sense
of unique principal asymptotic value for all $\zeta\in\partial D$
except a set of logarithmic capacity zero.

\section{On the logarithmic capacity}

The most important notion for our research is the notion of
logarithmic capacity, see e.g. \cite{No}, \cite{Kar} and \cite{N}.
First of all, given a bounded Borel set  $E$ in the plane $\Bbb C$,
a {\bf mass distribution} on $E$ is a nonnegative completely
additive function of a set $\nu$ defined on its Borel subsets with
$\nu(E)=1$. The function
\begin{equation}\label{4}
U^{\nu}(z):=\int\limits_{E} \log\left|\frac1{z-\zeta}
\right|\,d\nu(\zeta)
\end{equation}
is called a {\bf logarithmic potential} of the mass distribution
$\nu$ at a point $z\in\Bbb C$. A {\bf logarithmic capacity} $C(E)$
of the Borel set $E$ is the quantity
\begin{equation}\label{5}
C(E)=e^{-V}\ , \qquad V\ =\ \inf_{\nu}\ V_{\nu}(E)\ , \qquad
V_{\nu}(E)\ =\ \sup_{z}\ U^{\nu}(z)\ .
\end{equation}
Note that it is sufficient to take the supremum in (\ref{5}) over
the set $E$ only. If $V=\infty$, then $C(E)=0$. It is known that
$0\leq C(E)<\infty$, $C(E_1)\leq C(E_2)$ if $E_1\subseteq E_2$,
$C(E)=0$ if $E=\bigcup\limits_{n=1}\limits^{\infty } E_n$, with
$C(E_n)=0$, $n=1,2,\ldots$, see e.g. Lemma III.4 in \cite{Kar}.

It is also well-known the following geometric characterization of
the logarithmic capacity, see e.g. the point 110 in \cite{N}:
\begin{equation}\label{6}
C(E)\ =\ \tau(E)\ :=\ \lim_{n\to\infty}\ V_n^{\frac2{n(n-1)}}
\end{equation}
where $V_n$ denotes the supremum (really, maximum) of the product
\begin{equation}\label{7}
V(z_1,\dots,z_n)\ =\prod_{k<l}^{l=1,\dots,n}|z_k-z_l|
\end{equation}
taken over all collections of points $z_1,\dots,z_n$ in the set $E$.
Following F\'ekete, see \cite{F}, the quantity $\tau(E)$ is called
the {\bf transfinite diameter} of the set $E$. By the geometric
interpretation of the logarithmic capacity as the transfinite
diameter we immediately see that if $C(E)=0$, then $C(f(E))=0$ for
an arbitrary mapping $f$ that is continuous by H\"older and, in
particular, for conformal and quasiconformal mappings on the compact
sets, see e.g. Theorem II.4.3 in \cite{LV}.

In order to introduce sets that are measurable with respect to
logarithmic capacity, we define, following \cite{Kar}, {\bf inner
$C_*$ and outer $C^*$ capacities: }
\begin{equation}\label{INNER}
C_*(E)\ \colon =\ \sup_{F\subseteq E}\ C(E)
\end{equation}
where supremum is taken over all compact sets $F\subset\Bbb C$, and
\begin{equation}\label{OUTER}
C^{*}(E)\ \colon =\ \inf_{E\subseteq O}\ C(O)
\end{equation}
where infimum is taken over all open sets $O\subset\Bbb C$. Further,
a bounded set $E\subset\Bbb C$ is called {\bf measurable with
respect to the logarithmic capacity} if
\begin{equation}\label{MES}
C^{*}(E)\ =\ C_*(E)\ ,
\end{equation}
and the common value of $C_*(E)$ and $C^*(E)$ is still denoted by
$C(E)$. Note, see e.g. Lemma III.5 in \cite{Kar}, that the outer
capacity is semiadditive, i.e.,
\begin{equation}\label{ADD}
C^{*}\left(\bigcup\limits_{n=1}\limits^{\infty } E_n\right)\ \le\
\sum\limits_{n=1}\limits^{\infty } C^*(E_n)\ .
\end{equation}

A function $\varphi:E\to\mathbb{C}$ defined on a bounded set
$E\subset\mathbb{C}$ is called {\bf measurable with respect to
logarithmic capacity} if, for all open sets $O\subseteq\Bbb C$,  the
sets
\begin{equation}\label{FUN}
\Omega =\{ z\in E: \varphi(z)\in O\}
\end{equation}
are measurable with respect to logarithmic capacity. It is clear
from the definition that the set E is itself measurable with respect
to logarithmic capacity.

\medskip

Note also that sets of logarithmic capacity zero coincide with sets
of the so-called {\bf absolute harmonic measure} zero introduced by
Nevanlinna, see Chapter V in \cite{N}. Hence a set $E$ is of
(Hausdorff) length zero if $C(E)=0$, see Theorem V.6.2 in \cite{N}.
However, there exist sets of length zero having a positive
logarithmic capacity, see e.g. Theorem IV.5 in \cite{Kar}.

\medskip

\begin{remark}\label{R1}
It is known that Borel sets and, in particular, compact and open
sets are measurable with respect to logarithmic capacity, see e.g.
Lemma I.1 and Theorem III.7 in \cite{Kar}. Moreover, as it follows
from the definition, any set $E\subset\Bbb C$ of finite logarithmic
capacity can be represented as a union of the sigma-compactum (union
of countable collection of compact sets) and the set of logarithmic
capacity zero. It is also known that the Borel sets and, in
particular, compact sets are measurable with respect to all
Hausdorff's measures and, in particular, with respect to measure of
length, see e.g. theorem II(7.4) in \cite{S}. Consequently, any set
$E\subset\Bbb C$ of finite logarithmic capacity is measurable with
respect to measure of length. Thus, on such a set any function
$\varphi:E\to\mathbb{C}$ being measurable with respect to
logarithmic capacity is also measurable with respect to measure of
length on $E$. However, there exist functions that are measurable
with respect to measure of length but not measurable with respect to
logarithmic capacity, see e.g. Theorem IV.5 in \cite{Kar}.
\end{remark}

\medskip

We are especially interested by functions $\varphi:\partial\mathbb
D\to\mathbb{C}$ defined on the unit circle $\partial\mathbb D = \{
z\in\mathbb{C}: |z|=1\}$. However, in view of~\eqref{6}, it suffices
to examine the corresponding problems on the segments of the real
axis because any closed arc on $\partial \Bbb D$ admits a
bi-Lipschitz (even infinitely smooth, the so-called stereographic)
mapping $g$ onto such a segment.

\medskip

In this connection, recall that a mapping $g: X \to X'$ between
metric spaces $(X,d)$ and $(X',d')$ is said to be {\bf Lipschitz} if
$d'(g(x_1),g(x_2))\leqslant C\cdot d(x_1,x_2)$ for any $x_1,x_2 \in
X$ and for a finite constant $C$. If, in addition,
$d(x_1,x_2)\leqslant c\cdot d'(g(x_1),g(x_2))$ for any $x_1,x_2 \in
X$ and for a finite constant $c$, then mapping $g$ is called {\bf
bi-Lipschitz}.

\medskip

Recall also, see e.g. Subsection IV.10 in \cite{S}, that a point
$x_0\in\mathbb{R}$ is called a {\bf density point} for a measurable
(with respect to the length, i.e., with respect to the Lebesgue
measure) set $E\subset\mathbb{R}$ if $x_0\in E$ and
\begin{equation}\label{13}
\lim\limits_{\varepsilon\to0}\ \frac{|\ (x_0-\varepsilon,
x_0+\varepsilon)\setminus E\ |}{2\varepsilon}\ =\ 0\ .
\end{equation}

\medskip

Similarly, we say that a point $x_0\in\mathbb{R}$ is the {\bf
density point with respect to logarithmic capacity} for a measurable
(with respect to $C$) set $E\subset\mathbb{R}$ if $x_0\in E$ and
\begin{equation}\label{D}
\lim\limits_{\varepsilon\to0}\ \frac{C([x_0-\varepsilon,
x_0+\varepsilon]\setminus E)}{C([x_0-\varepsilon,
x_0+\varepsilon])}\ =\ 0\ .
\end{equation}

\medskip

Note that the logarithmic capacity of a segment of length $l$ is
equal to $l/4$, see e.g. \cite{ST}, p. 25 and 45. Hence
$C([x_0-\varepsilon, x_0+\varepsilon])=\varepsilon / 2$. Note
simultaneously from the same place that logarithmic capacity of a
circle (and a disk) of a radius $r$ is equal to $r$ and, in
particular, logarithmic capacity of the unit circle (and the unit
disk) is equal to 1.

\medskip

Finally, recall that a function $\varphi:[a,b]\to\mathbb{C}$ is {\bf
approximately continuous (with respect to logarithmic capacity) at a
point} $x_0\in(a,b)$ if it is continuous on a set $E\subseteq[a,b]$
for which $x_0$ is a density point (with respect to logarithmic
capacity), see e.g. Subsection IV.10 in \cite{S}, and the Subsection
2.9.12 in \cite{Fe}, correspondingly.

\medskip

Further, it is important that the following analog of the Denjoy
theorem holds, see e.g. Theorem 2.9.13 in \cite{Fe}, cf. Theorem
IV(10.6) in \cite{S}.

\medskip

\begin{proposition}\label{P1} {\it
A function $\varphi:[a,b]\to\mathbb{C}$ is measurable with respect
to logarithmic capacity if and only if it is approximately
continuous for a.e. $x\in(a,b)$ with respect to logarithmic
capacity.}
\end{proposition}

\medskip

\begin{remark}\label{R2}
As known, $C([a; b]) = (b - a)/4$ and, moreover, $C(E) \ge |E|/4 $
where $|E|$ – длина $E$, see e.g. II.4.17 in \cite{La}. Thus, if
$x_0$ is a density point for a set $E$ with respect to logarithmic
capacity, then $x_0$ is also a density point for the set $E$ with
respect to measure of length. Consequently, each point of
approximate continuity of a function $\varphi:[a,b]\to\Bbb R$ with
respect to logarithmic capacity is also the point of approximate
continuity of the function $\varphi$ with respect to the Lebesgue
measure on the real axis.\end{remark}

\medskip

Hence, in particular, we obtain the following useful lemma.

\medskip

\begin{lemma}\label{L1}
{\it Let a function $\varphi:[a,b]\to\Bbb R$ be bounded and
measurable with respect to logarithmic capacity and let
$\Phi(x)=\int\limits_a\limits^x\varphi(t)\ dt$ be its indefinite
Lebesgue integral. Then $\Phi^{\prime}(x)=\varphi(x)$ a.e. on
$(a,b)$ with respect to logarithmic capacity.}
\end{lemma}

\medskip

\begin{proof}
Indeed, let $x_0\in(a,b)$ be a point of approximate continuity for a
function $\varphi$. Then there is a set $E\subseteq[a,b]$ for which
point $x_0$ is a density point and $\varphi$ is continuous on this
set. Since $|\varphi(x)|\le C<\infty$ for all $x\in[a,b]$, we obtain
that at small $h$
$$
\left | \frac{\Phi(x_0+h)-\Phi(x_0)}{h}-\varphi(x_0)\right |\ \le\
$$
$$
\le\ \max\limits_{x\in E\cap[x_0,x_0+h]}\ |\
\varphi(x)-\varphi(x_0)\ |\ + 2C\ \frac{|\ (x_0,x_0+h)\setminus E\
|}{|h|}\ ,
$$
i.e., $\Phi^{\prime}(x_0)=\varphi(x_0)$. Thus, Lemma \ref{L1}
follows from Proposition \ref{P1}, see also Remarks \ref{R1} and
\ref{R2}.
\end{proof}

\section{One analog of the Lusin theorem}

The following remarkable theorem of Lusin says that, for any
measurable finite  a.e. (with respect to the Lebesgue measure)
function $\varphi$ on the segment $[a,b]$, there is a continuous
function $\Phi$ such that $\Phi^{\prime}(x)=\varphi(x)$ a.e. on
$[a,b]$, see e.g. Theorem VII(2.3) in \cite{S}. This statement was
well-known long ago for integrable functions $\varphi$ with respect
to its indefinite integral $\Phi$, see e.g. Theorem IV(6.3) in
\cite{S}. However, this result is completely nontrivial for
nonintegrable functions $\varphi$.

In the proof of one analog of the Lusin theorem in terms of
logarithmic capacity, the following lemma on singular functions of
the Cantor type will take the key part.

\medskip

\begin{lemma}\label{L2} {\it
There is a continuous nondecreasing function $\Psi:[0,1]\to[0,1]$
such that $\Psi(0)=0$, $\Psi(1)=1$ and $\Psi^{\prime}(t)=0$ a.e.
with respect to logarithmic capacity.}
\end{lemma}

\medskip

\begin{proof}
To prove this fact we use the construction of sets of Cantor's type
of logarithmic capacity zero due to Nevanlinna. Namely, let us
consider a sequence of numbers $p_k>1$, $k=1, 2, \ldots$, and define
the corresponding sequence of the sets $E(p_1,\ldots, p_n)$, $n=1,
2, \ldots$, by the induction in the following way. Let $E(p_1)$ be
the set consisting of two equal-length segments obtained from the
unit segment $[0,1]$ by removing the central interval of length
$1-1/p_1$; $E(p_1, p_2)$ be the set consisting of $2^2=4$
equal-length segments obtained by removing from each segment of the
previous set $E(p_1)$ the central interval with $1-1/p_2$ fraction
of its length and so on. Denote by $E(p_1, p_2, \ldots)$ the
intersection of all the sets $E(p_1,\ldots, p_n)$, $n=1, 2, \ldots$.
By Theorem V.6.3 in \cite{N} the set $E(p_1, p_2, \ldots)$ has
logarithmic capacity zero if and only if the series $\sum 2^{-k}\log
p_k$ is divergent. This condition holds, for example, if $p_k=e^{
2^k }$.

It is known that all sets of Cantor's type are homeomorphic each to
other. In particular, there is a homeomorphism $h:[0,1]\to[0,1]$,
$h(0)=0$ and $h(1)=1$, under which $E(p_1, p_2, \ldots)$ is
transformed into the classical Cantor set, see e.g. 8.23 in
\cite{GO}. Thus, if $\varkappa$ is the classical Cantor function,
see e.g. 8.15 in \cite{GO}, then $\Psi=\varkappa\circ h$ is the
desired function.
\end{proof}$\Box$

\medskip

\begin{lemma}\label{L3} {\it
Let a function $g:[a,b]\to\mathbb{R}$ be bounded and measurable with
respect to logarithmic capacity. Then, for every $\varepsilon>0$,
there is a continuous function $G : [a,b]\to\mathbb{R}$ such that
$|G(x)|\le\varepsilon$ for all $x\in[a,b]$,  $G(a)=G(b)=0$, and
$G^{\prime}(x)=g(x)$ a.e. on $[a,b]$ with respect to logarithmic
capacity.}
\end{lemma}

\medskip

\begin{proof}
Let $H(x)=\int\limits_a\limits^xg(t)\ dt$ be the indefinite Lebesgue
integral of the function $g$. Choose on $[a,b]$ a finite collection
of points $a=a_0<a_1<\ldots < a_n=b$ such that the oscillation of
$H$ on each segment $[a_k, a_{k+1}]$, $k=0,1,\ldots , n-1$ is less
than $\varepsilon /2$. Applying linear transformations of
independent and dependent variables to the function
$\Psi:[0,1]\to[0,1]$ from Lemma \ref{L2}, we obtain the function
$F_k$ on each segment $[a_k, a_{k+1}]$, $k=0,1,\ldots , n-1$, that
coincides with the function $H$ at its endpoints and whose
derivative is equal to zero a.e. with respect to logarithmic
capacity. Let $F$ be the function on $[a,b]$ glued of the functions
$F_k$. Then $G=H-F$ gives us the desired function by Lemmas \ref{L1}
and \ref{L2}.
\end{proof}$\Box$

\medskip

\begin{lemma}\label{L4} {\it
Let a function $g:[a,b]\to\mathbb{R}$ be bounded and measurable with
respect to logarithmic capacity and let $P$ be a closed subset of
the segment $[a,b]$. Then, for every $\varepsilon>0$, there is a
continuous function $G : [a,b]\to\mathbb{R}$ such that
$|G(x+h)|\le\varepsilon |h|$ for all $x\in P$ and all
$h\in\mathbb{R}$ such that $x+h\in[a,b]$, $G(x)=G^{\prime}(x)=0$ for
all $x\in P$ and $G^{\prime}(x)=g(x)$ a.e. on $[a,b]\setminus P$
with respect to logarithmic capacity.}
\end{lemma}

\medskip

\begin{proof}
Let $I=(a,b)$. Then the set $I\setminus P$ is open and can be
represented as the union of a countable collection of mutually
disjoint intervals $I_k=(a_k,b_k)$. Choose in each interval $I_k$ an
increasing sequence of numbers $с^{(j)}_k$, $j=0, \pm 1, \pm 2,
\ldots$ such that $c^{(j)}_k\to a_k$  as $j\to -\infty$ and
$c^{(j)}_k\to b_k$ as $j\to +\infty$. Denote by
$\varepsilon^{(j)}_k$ the minimal of the two numbers $\varepsilon
(c^{(j)}_k-a_k)/(k+|j|)$ and $\varepsilon (b_k-c^{(j)}_k)/(k+|j|)$.
Then by Lemma \ref{L3} in each interval $I_k$ there is a continuous
function $G_k$ such that $|G(x)|\le\varepsilon^{(j)}_k$ for all
$x\in[c^{(j)}_k, c^{(j+1)}_k]$, $G(c^{(j)}_k)=0$ for all $j=0, \pm
1, \pm 2, \ldots$, and $G^{\prime}(x)=g(x)$ a.e. on $I_k$ with
respect to logarithmic capacity. Thus, setting $G(x)=G_k(x)$ on each
interval $I_k$ and $G(x)=0$ on the set $P$, we obtain the desired
function.
\end{proof}$\Box$

\medskip

Finally, we prove the following analog of the Lusin theorem
mentioned above.

\medskip

\begin{theorem}\label{T1} {\it
Let $\varphi:[a,b]\to\mathbb R$ be a measurable function with
respect to logarithmic capacity. Then there is a continuous function
$\Phi:[a,b]\to\mathbb R$ such that $\Phi^{\prime}(x)=\varphi(x)$
a.e. on $(a,b)$ with respect to logarithmic capacity. Furthermore,
the function $\Phi$ can be choosen such that $\Phi(a)=\Phi(b)=0$ and
$|\Phi(x)|\le\varepsilon$ for a prescribed $\varepsilon>0$ and all
$x\in[a,b]$.}
\end{theorem}

\medskip

\begin{proof}
First we define by induction a sequence of closed sets $P_n\subseteq
[a,b]$ and a sequence of continuous functions $G_n:[a,b]\to\mathbb
R$, $n=0, 1,\dots$, whose derivatives exist a.e. and are measurable
with respect to logarithmic capacity such that, under the notations
$Q_n=\bigcup\limits_{k=0}\limits^{n}P_k$ and
$\Phi_n=\sum\limits_{k=0}\limits^{n}G_k$, the following conditions
hold: $ (a)\ \Phi_n^{\prime}(x) = \varphi(x) \ \mbox{for}\ x\in
Q_n,\ $ $ (b)\  G_n(x) = 0 \ \mbox{for}\ x\in Q_{n-1},\ $ $ (c)\
|G_n(x+h)|\le |h|/2^n $ for all $ x\in Q_{n-1}\ \mbox{and}$ all $h$
such that $ x+h\in[a,b],\ $ $(d)\ C(I\setminus Q_n)<1/n$ where
$I=[a,b]$.

So, let $G_0\equiv 0$ and $P_0=\varnothing$ and let $G_n$ and $P_n$
be already constructed with the given conditions for all $n=1,
2,\dots, m$. Then there is a compact set $E_m\subset I\setminus Q_m$
such that
\begin{equation}\label{15}
 C(I\setminus(Q_m\cup
E_m))<1/(m+1)
\end{equation}
and the functions $\Phi_m^{\prime}$ and $\varphi$ are continuous on
$E_m$, see e.g. Theorem 2.3.5 in \cite{Fe}.

By Lemma \ref{L4} with the set $P=Q_m$ and the function
$g:I\to\mathbb{R}$ that is equal to $\varphi(x)-\Phi_m^{\prime}(x)$
on $E_m$ and zero on $I\setminus E_m$, there is a continuous
function $G_{m+1}: I\to\mathbb{R}$ such that $(i)\
G_{m+1}^{\prime}(x)=\varphi(x)-\Phi_m^{\prime}(x)$ a.e. on
$I\setminus Q_m$ with respect to logarithmic capacity, $(ii)\
G_{m+1}(x)=G_{m+1}^{\prime}(x)=0$ for all $x\in Q_m$, and $(iii)\
|G_{m+1}(x+h)|\le |h|/2^{m+1}$ for all $x\in Q_m$ and all $h$ such
that $x+h\in I$.

By the definition of logarithmic capacity, conditions $(i)$ and
(\ref{15}), there is a compact set $P_{m+1}\subseteq E_m$ such that
\begin{equation}\label{16}
C(I\setminus(Q_m\cup P_{m+1}))<1/(m+1)\ ,
\end{equation}
\begin{equation}\label{17}
G_{m+1}^{\prime}(x)=\varphi(x)-\Phi_m^{\prime}(x)\ \ \ \ \ \forall\
x\in P_{m+1} \ .
\end{equation}
By (\ref{16}), (\ref{17}), $(ii)$, $(iii)$ it is easy to see that
conditions $(a)$, $(b)$, $(c)$ and $(d)$ hold for $n=m+1$, too.

Set now on the basis of the above construction of the sequences
$G_n$ and $P_n$:
\begin{equation}\label{18}
\Phi(x)\ =\ \lim\limits_{k\to\infty}\ \Phi_k(x)\ =\
\sum\limits_{k=1}\limits^{\infty}\ G_k(x)\ , \ \ \ \ \ \ Q\ =\
\lim\limits_{k\to\infty}\ Q_k\ =\
\bigcup\limits_{k=1}\limits^{\infty}\ P_k\ .
\end{equation}
Note that $\Phi_k\to\Phi$ uniformly on the segment $I$ because of
the condition $(c)$ and hence the function $\Phi$ is continuous. By
the construction, for each $x_0\in Q$ we have that $x_0\in Q_n$ for
large enough $n$ and, since
$$
\frac{\Phi(x_0+h)-\Phi(x_0)}{h}\ =\
\frac{\Phi_n(x_0+h)-\Phi_n(x_0)}{h}\ +\
\sum\limits_{k=n+1}\limits^{\infty}\ \frac{G_k(x_0+h)-G_k(x_0)}{h}\
,
$$
we obtain from the conditions $(a)$, $(b)$ and $(c)$ that
$$
\limsup\limits_{h\to0}\ \left |\ \frac{\Phi(x_0+h)-\Phi(x_0)}{h}\ -
\varphi(x_0)\ \right |\ <\ \frac{1}{2^n}\ ,
$$
i.e., $\Phi^{\prime}(x_0)=\varphi(x_0)$. Moreover, by condition
$(d)$ we see that $C(I\setminus Q)=0$. Thus,
$\Phi^{\prime}(x)=\varphi(x)$ a.e. on $[a,b]$ with respect to
logarithmic capacity.

Finally, applying the construction of the proof of Lemma \ref{L3} to
the function $\Phi$ instead of the indefinite integral, we find a
new function $\Phi_*$ such that $\Phi_*^{\prime}(x)=\varphi(x)$ a.e.
on $[a,b]$ with respect to logarithmic capacity with
$\Phi_*(a)=\Phi_*(b)=0$ and $|\Phi_*(x)|\le\varepsilon$ for a
prescribed $\varepsilon>0$ and all $x\in[a,b]$.
\end{proof}$\Box$

\section{On the Dirichlet problem for harmonic functions in the unit circle}

Gehring in \cite{Ge} has established the following brilliant result:
if $\varphi:\Bbb R\to\Bbb R$ is $2\pi$-periodic, measurable and
finite a.e. with respect to the Lebesgue measure, then there is a
harmonic function in $|z|<1$ such that $u(z)\to\varphi(\vartheta)$
for a.e. $\vartheta$ as $z\to e^{i\vartheta}$ along any
nontangential path.

\medskip

It will be useful the following analog of the Gehring theorem.

\medskip

\begin{theorem}\label{T2} {\it
Let $\varphi:\Bbb R\to\Bbb R$ be $2\pi$-periodic, measurable and
finite a.e. with respect to logarithmic capacity. Then there is a
harmonic function $u(z)$, $z\in\Bbb D$, such that
$u(z)\to\varphi(\vartheta)$ as $z\to e^{i\vartheta}$ along any
nontangential path for all $\vartheta\in\Bbb R$ except a set of
logarithmic capacity zero.}
\end{theorem}

\medskip

\begin{proof}
By Theorem \ref{T1} we are able to find a continuous $2\pi$-periodic
function $\Phi:\Bbb R\to\Bbb R$ such that
$\Phi^{\prime}(\vartheta)=\varphi(\vartheta)$ for a.e. $\vartheta$
with respect to logarithmic capacity. Set
\begin{equation}\label{Poisson}
U(re^{i\vartheta})\ =\ \frac{1}{2\pi}\
\int\limits_{0}\limits^{2\pi}\frac{1-r^2}{1-2r\cos(\vartheta-t)+r^2}\
\Phi(t)\ dt
\end{equation}
for $r<1$. Next, by the well-known result due to Fatou, see e.g.
3.441 in \cite{Z}, p. 53, see also Theorem IX.1.2 in \cite{Go},
$\frac{\partial}{\partial\vartheta}\ U(z)\to
\Phi^{\prime}(\vartheta)$ as $z\to e^{i\vartheta}$ along any
nontangential path whenever $\Phi^{\prime}(\vartheta)$ exists. Thus,
the conclusion follows for the function $u(z)\ =\
\frac{\partial}{\partial\vartheta}\ U(z)$.
\end{proof}$\Box$

\medskip

\begin{remark}\label{remZERO} Note that the given function $u$ is
harmonic in the punctured unit disk $\mathbb D\setminus\{0\}$
because the function $U$ is harmonic in $\mathbb D$ and the
differential operator $\frac{\partial}{\partial\vartheta}$ is
commutative with the Laplace operator $\Delta$. Setting $u(0)=0$, we
see that
$$
u(re^{i\vartheta})\ =\ -\frac{r}{\pi}\
\int\limits_{0}\limits^{2\pi}\frac{(1-r^2)\sin(\vartheta-t)}{(1-2r\cos(\vartheta-t)+r^2)^2}\
\Phi(t)\ dt\ \to\ 0\ \ \ \ \ \ \ \mbox{as}\ r\to 0\ ,
$$
i.e. $u(z)\to u(0)$ as $z\to 0$, and, moreover, the integral of $u$
over each circle $|z|=r$, $0<r<1$, is equal to zero. Thus, by the
criterion for a harmonic function on the averages over circles we
have that $u$  is harmonic in $\mathbb D$. The alternative argument
for the latter is the removability of isolated singularities for
harmonic functions, see e.g. \cite{N}.\end{remark}

\medskip

It is known that every harmonic function $u(z)$ in $\mathbb D = \{
z\in\mathbb C : |z|<1\}$ has a conjugate function $v(z)$ such that
$f(z)=u(z)+iv(z)$ is an analytic function in $\mathbb D$. Hence we
have the following corollary:

\medskip

\begin{corollary}\label{C1} {\it
Under the conditions of Theorem \ref{T2}, there is an analytic
function $f$ in $\mathbb D$ such that $\rm Re\ f(z)\to
\varphi(\vartheta)$ as $z\to e^{i\vartheta}$ along any nontangential
path for a.e. $\vartheta$ with respect to logarithmic capacity.}
\end{corollary}

\medskip

Note that the boundary values of the conjugate function $v$ cannot
be prescribed arbitrarily and simultaneously with the boundary
values of $u$ because $v$ is uniquely determined by $u$ up to an
additive constant, see e.g. I.A in \cite{Ku}.

\medskip

Denote by $h^p$, $p\in(0,\infty)$, the class of all harmonic
functions $u$ in $\mathbb D$ with
$$
\sup\limits_{r\in(0,1)}\ \left\{\int\limits_{0}\limits^{2\pi}\
|u(re^{i\vartheta})|^p\ d\vartheta\right\}^{\frac{1}{p}}\ <\ \infty\
.
$$

\medskip

\begin{remark}\label{R3}
It is clear that $h^p\subseteq h^{p^{\prime}}$ for all
$p>p^{\prime}$ and, in particular, $h^p\subseteq h^{1}$ for all
$p>1$. It is important that every function in the class $h^1$ has
a.e. nontangential boundary limits, see e.g. Corollary IX.2.2 in
\cite{Go}.
\end{remark}

\medskip

It is also known that a harmonic function $U$ in $\Bbb D$ can be
represented as the Poisson integral~\eqref{Poisson} with a function
$\Phi\in L^p(-\pi,\pi),\: p>1$, if and only if $U\in h^p$, see e.g.
Theorem IX.2.3 in \cite{Go}. Thus, $U(z)\to \Phi(\vartheta)$ as
$z\to e^{i\vartheta}$ along any nontangential path for a.e.
$\vartheta$, see e.g. Corollary IX.1.1 in \cite{Go}. Moreover,
$U(z)\to \Phi(\vartheta_0)$ as $z\to e^{i\vartheta_0}$ at points
$\vartheta_0$ of continuity of the function $\Phi$, see e.g. Theorem
IX.1.1 in \cite{Go}.

Note also that $v\in h^p$ whenever $u\in h^p$ for all $p>1$ by the
M. Riesz theorem, see \cite{R}. Generally speaking, this fact is not
trivial but it follows immediately for $p=2$ from the Parseval
equality, see e.g. the proof of Theorem IX.2.4 in \cite{Go}. The
latter will be sufficient for our goals.

\section{Correlations of boundary data of conjugate functions}

It is known the very delicate fact due to Lusin that harmonic
functions in the unit circle with continuous (even absolutely
continuous !) boundary data can have conjugate harmonic functions
whose boundary data are not continuous functions, furthemore, they
can be even not essentially bounded in neighborhoods of each point
of the unit circle, see e.g. Theorem VIII.13.1 in \cite{Bari}. Thus,
a correlation between boundary data of conjugate harmonic functions
is not a simple matter, see also I.E in \cite{Ku}.

\medskip

We call $\lambda:\partial\mathbb D\to\mathbb C$ a {\bf function of
bounded variation}, write $\lambda\in\mathcal{BV}(\partial\mathbb
D)$, if
\begin{equation}\label{20}
V_{\lambda}(\partial\mathbb D)\ \colon =\ \sup\
\sum\limits_{j=1}\limits^{j=k}\
|\lambda(\zeta_{j+1})-\lambda(\zeta_j)| \ <\ \infty \end{equation}
where the supremum is taken over all finite collections of points
$\zeta_j\in\partial\mathbb D$, $j=1,\ldots , k$, with the cyclic
order meaning that $\zeta_j$ lies between $\zeta_{j+1}$ and
$\zeta_{j-1}$ for every $j=1,\ldots , k$. Here we assume that
$\zeta_{k+1}=\zeta_1=\zeta_0$. The quantity
$V_{\lambda}(\partial\mathbb D)$ is called the {\bf variation of the
function} $\lambda$.

\medskip

\begin{remark}\label{R4}
It is clear by the triangle inequality that if we add new
intermediate points in the collection $\zeta_j$, $j=1,\ldots , k$,
then the sum in \eqref{20} does not decrease. Thus, the given
supremum is attained as $\delta=\sup\limits_{j=1,\ldots
k}|\zeta_{j+1}-\zeta_j|\to 0$. Note also that by the definition
$V_{\lambda}(\partial\mathbb D)=V_{\lambda\circ h}(\partial\mathbb
D)$, i.e., the {\bf variation is invariant} under every
homeomorphism $h:\partial\mathbb D\to\partial\mathbb D$ and, thus,
the definition can be extended in a natural way to an arbitrary
Jordan curve in $\mathbb C$.
\end{remark}

\medskip

\begin{lemma}\label{T3} {\it
Let $\mathbb D$ be the unit disk in $\mathbb C$,
$\alpha:\partial\mathbb D\to\mathbb R$ be a function of bounded
variation, $u:\mathbb D\to\mathbb R$ be a bounded harmonic function
such that
\begin{equation}\label{23}
\lim\limits_{z\to\zeta}\ u(z)\ =\ \alpha(\zeta)
\end{equation}
at every point of continuity of $\alpha$ and let $v$ be its
conjugate harmonic function. Then for a.e. $\zeta\in\partial\mathbb
D$ with respect to logarithmic capacity
\begin{equation}\label{24}
\lim\limits_{z\to\zeta}\ v(z)\ =\ \beta(\zeta)
\end{equation}
along any nontangential path in $\mathbb D$ terminating at $\zeta$
where $\beta:\partial\mathbb D\to\mathbb R$ is a function that is
measurable with respect to logarithmic capacity.}
\end{lemma}

\medskip

\begin{proof}
Indeed,  $\alpha$ as a function of bounded variation has at most a
countable set $S$ of points of discontinuity and, consequently,
$C(S)=0$. Hence by the generalized maximum principle, see e.g. the
point 115 in \cite{N}, such a function $u$ is unique. Moreover,
$\alpha\in\mathcal{BV}(\partial\mathbb D)$ is bounded and by the
Denjoy theorem, see e.g. Theorem IV(10.6) in \cite{S}, cf.
Proposition \ref{P1} above, it is measurable with respect to the
length measure (as well as with respect to logarithmic capacity),
i.e., $\alpha\in L^{\infty}(\partial\mathbb D)$, and, consequently,
$u$ can be represented as the Poisson integral of the function
$\alpha$, see e.g. Theorem I.D.2.2 in \cite{Ku},
\begin{equation}\label{PoissonV}
u(re^{i\vartheta})\ =\ \frac{1}{2\pi}\
\int\limits_{-\pi}\limits^{\pi}\frac{1-r^2}{1-2r\cos(\vartheta
-t)+r^2}\ \alpha(e^{it})\ dt\ .
\end{equation}
Here the Poisson kernel is a real part of the analytic function
$(\zeta +z)/(\zeta -z)$, $\zeta = e^{it}$, $z=re^{i\vartheta}$, and
by the Weierstrass theorem, see e.g. Theorem 1.1.1 in \cite{Go}, the
Schwartz integral
\begin{equation}\label{PoissonW}
f(z)\ :=\ \frac{1}{2\pi i}\ \int\limits_{\partial\mathbb D}
\alpha(\zeta)\ \frac{\zeta +z}{\zeta -z}\ \frac{d\zeta}{\zeta}
\end{equation}
gives the analytic function $f=u+iv$ in $\mathbb D$ with  $u={\rm
Re}\, f$, $v={\rm Im}\, f$, and
\begin{equation}\label{PoissonT}
f(z)\ =\ \frac{1}{2\pi}\ \int\limits_{-\pi}\limits^{\pi}
\alpha(e^{it})\ \frac{e^{it} +z}{e^{it} -z}\ dt\ =\ C\ +\
\frac{z}{\pi}\ \int\limits_{-\pi}\limits^{\pi}
\frac{F(t)}{1-e^{-it}z}\ dt
\end{equation}
where $F(t)=e^{-it}\alpha(e^{it})$ and $C=\frac{1}{2\pi}
\int\limits_{-\pi}\limits^{\pi} \alpha(e^{it})\ dt$. By Theorem 2(c)
in \cite{T} the function $f(z)$ has angular limits  $f(\zeta)$ as
$z\to\zeta$ for a.e. $\zeta\in\partial\mathbb D$ with respect to
logarithmic capacity because the function $F$ is of bounded
variation. It remains to note that
$f(\zeta)=\lim\limits_{n\to\infty}\, f_n(\zeta)$, where
$f_n(\zeta)=f(r_n\zeta)$, for an arbitrary sequence $r_n\to 1-0$ as
$n\to \infty$ for a.e. $\zeta\in\partial\mathbb D$ with respect to
logarithmic capacity and, thus, $f(\zeta)$ is measurable with
respect to logarithmic capacity because the functions $f_n(\zeta)$
are so as continuous functions on $\partial\mathbb D$, see e.g.
2.3.10 in \cite{Fe}.
\end{proof}$\Box$

\medskip

\begin{remark}\label{RT}
One can show that the statement is valid for analytic functions in
quasidisks $D$ and, in particular, in Jordan domains with smooth and
Lipschitz boundaries. By the given proof, see especially Theorem
2(c) in \cite{T}, it is also clear that $f=u+iv$ admits limits along
wide classes of tangential paths to a.e. boundary point of
$\partial\mathbb D$ with respect to logarithmic capacity. However,
it is not interesting for us in the context of our research.
\end{remark}

\medskip

We also prove the following statement that will be useful later on.

\medskip

\begin{proposition}\label{P2} {\it
For every function $\lambda:\partial\mathbb{D}\to\partial\mathbb{D}$
of the class $\mathcal{BV}(\partial\mathbb{D})$ there is a function
$\alpha_{\lambda}:\partial\mathbb{D}\to\mathbb{R}$ of the class
$\mathcal{BV}(\partial\mathbb{D})$ with $V_{\alpha_{\lambda}}\le
V_{\lambda}\cdot 3\pi/2$ such that
$\lambda(\zeta)=\exp\{{i\alpha_{\lambda}}(\zeta)\}$,
$\zeta\in\partial\mathbb{D}$.}
\end{proposition}

\medskip

We will call the function $\alpha_{\lambda}$ a {\bf function of
argument} of $\lambda$.

\medskip

\begin{proof}
Let us consider the function
$\Lambda(\vartheta)=\lambda(e^{i\vartheta})$,
$\vartheta\in[0,2\pi]$. It is clear that $V_{\Lambda}=V_{\lambda}$
and, thus, $\Lambda$ has not more than a countable collection of
jumps $j_n$ where the series $\sum j_n$ is absolutely convergent,
$\sum |j_n|\le V_{\lambda}$, and
$\Lambda(\vartheta)=J(\vartheta)+C(\vartheta)$ where $C(\vartheta)$
is a continuous function and $J(\vartheta)$ is the function of jumps
of $\Lambda$ that is equal to the sum of all its jumps in
$[0,\vartheta]$, see e.g. Corollary VIII.3.2 and Theorem VIII.3.7 in
\cite{Nat}. We have that $V_J\le V_{\lambda}$ and $V_C\le
2V_{\lambda}$, see e.g. Theorem 6.4 in \cite{Rud}. Let us associate
with the complex quantity $j_n$ the real quantity
$$\alpha_n\ =\ -2\ {\rm arctg}\ \frac{{\rm Re}\ j_n}{{\rm Im}\ j_n}\ \in\ [-\pi,\pi]\ .$$ By the
geometric interpretation of these quantities ($|j_n|$ is equal to
the length of the chord for an arc of the unit circle of the length
$|\alpha_n|$) and elementary calculations, we have that $|j_n| \le
|\alpha_n| \le |j_n|\cdot\pi/2$.

The first inequality is evident because the length $|\alpha|$ of an
arc $\alpha$ of a circle is always greater than the length $|h|$ of
a chord $h$ connecting its ends. For the proof of the second
inequality, note that $|h|=2|\sin \frac{\alpha}{2}|$. The necessary
condition for extremum of the functional $f(\beta):=\frac{\sin
^2\beta}{\beta ^2}$, $\beta\in[0,\frac{\pi}{2}]$ in the points of
its differentiability is the equality
$f^{\prime}(\beta_0)=0=\frac{2\sin\beta_0}{\beta_0^3}(\beta_0\cos\beta_0
- \sin\beta_0)$, i.e. $\tan\beta_0=\beta_0$. However, the functional
$g(\beta):=\tan\beta -\beta$, $\beta\in[0,\frac{\pi}{2}]$, is
strictly increasing because $g^{\prime}(\beta)=\tan^2\beta>0$ for
$\beta\in(0,\frac{\pi}{2})$, $g(0)=0$, and hence the equality
$\tan\beta_0=\beta_0$ cannot hold for $\beta\in(0,\frac{\pi}{2})$.
Thus, the extreme points of the functional $f$ are the ends of the
interval.

Now, let us associate with the function $J(\vartheta)$ the function
$j(\vartheta)$ that is equal to the sum of all $\alpha_n$
corresponding to jumps of $\Lambda$ in $[0,\vartheta]$, $V_j\le
V_J\cdot\pi/2$. Next, let us associate with the complex-valued
function $C(\vartheta)$ a real-valued function $c(\vartheta)$ in the
following way. As $C(\vartheta)$ is uniformly continuous on the
segment $[0,2\pi]$, the latter can be split to the segments
$S_k=[\theta_{k-1},\theta_k]$, $\theta_k = 2\pi k/m$, $k=1,\ldots ,
m$, with a large enough $m\in\Bbb N$ such that
$|\Lambda(\vartheta)-\Lambda(\vartheta^{\prime})|<2$ for all
$\vartheta$ and $\vartheta^{\prime}\in S_{k}$. Set by induction
$$c(\vartheta)=c(\theta_{k-1})-2\ {\rm arctg}\ \frac{{\rm Re}
[C(\vartheta)-C(\theta_{k-1})]}{{\rm Im}
[C(\vartheta)-C(\theta_{k-1})]} \ \ \ \ \ \ \ \ \  \
\forall\vartheta\in S_k,\: k=1,\dots,m,$$ where
$$c(0):={\rm arctg}\ \frac{{\rm Re}[C(0)-1]}{{\rm Im}[C(0)-1]}\:\:.$$ Moreover, let
$\gamma_{\lambda}(\vartheta)=j(\vartheta)+c(\vartheta)$,
$\vartheta\in[0,2\pi]$. By the construction
$\Lambda(\vartheta)=e^{i\gamma_{\lambda}(\vartheta)}$,
$\vartheta\in[0,2\pi]$, $V_{\gamma_{\lambda}}\le V_{{\lambda}}\cdot
3\pi/2$. Finally, setting
$\alpha_{\lambda}(\zeta)=\gamma_{\lambda}(\vartheta)$ for
$\zeta=e^{i\vartheta}$, $\vartheta\in[0,2\pi)$, we obtain the
desired function $\alpha_{\lambda}$ of the class
$\mathcal{BV}(\partial\mathbb{D})$.
\end{proof}$\Box$

\bigskip

\section{The Riemann-Hilbert problem for analytic functions}

\begin{theorem}\label{T4} {\it
Let $\lambda:\partial\mathbb D\to\partial\mathbb D$ be of bounded
variation and $\varphi:\partial\mathbb D\to\mathbb R$ be measurable
with respect to logarithmic capacity. Then there is an analytic
function $f:\mathbb D\to\mathbb C$ such that along any nontangential
path
\begin{equation}\label{25}
\lim\limits_{z\to\zeta}\ \mathrm {Re}\
\{\overline{\lambda(\zeta)}\cdot f(z)\}\ =\ \varphi(\zeta)
\quad\quad\quad\mbox{for}\ \ \mbox{a.e.}\ \ \
\zeta\in\partial\mathbb D
\end{equation}
with respect to logarithmic capacity.}
\end{theorem}

\medskip

\begin{proof}
By Proposition \ref{P2} the function of argument
$\alpha_{\lambda}\in\mathcal{BV}(\partial\mathbb{D})$. Therefore
$$ g(z)\ =\ \frac{1}{2\pi i}\ \int\limits_{\partial\mathbb
D}\alpha(\zeta)\ \frac{z+\zeta}{z-\zeta}\
 \frac{d\zeta}{\zeta}\ , \ \ \ \ \ z\in\mathbb D\ ,
$$
is an analytic function with $u(z)={\mathrm Re}\
g(z)\to\alpha(\zeta)$ as $z\to\zeta$ for every
$\zeta\in\partial\mathbb D$ except a countable collection of points
of discontinuity of $\alpha_{\lambda}$, see e.g. Corollary IX.1.1 in
\cite{Go} and Theorem I.D.2.2 in \cite{Ku}. Note that ${\cal
A}(z)=\exp\{ig(z)\}$ is also an analytic function.

By Lemma \ref{T3} there is a function $\beta:\partial\mathbb
D\to\mathbb R$ that is finite a.e. and measurable with respect to
logarithmic capacity such that $v(z)={\mathrm Im}\
g(z)\to\beta(\zeta)$ as $z\to\zeta$ for a.e.
$\zeta\in\partial\mathbb D$ with respect to logarithmic capacity
along any nontangential path. Thus, by Corollary \ref{C1} there is
an analytic function ${\cal B}:\mathbb D\to\mathbb C$ such that
$U(z)=\mathrm {Re}\ {\cal B}(z)\to\varphi(\zeta)\cdot
\exp\{{\beta(\zeta)}\}$ as $z\to\zeta$ along any nontangential path
for a.e. $\zeta\in\partial\mathbb D$. Finally, elementary
calculations show that the desired function $f={\cal A}\cdot{\cal
B}$.
\end{proof}$\Box$

\medskip

By the Bagemihl theorem, see Introduction, we obtain directly from
Theorem \ref{T4} the following result.

\medskip

\begin{theorem}\label{T5} {\it
Let $D$ be a Jordan domain in $\Bbb C$, $\lambda:\partial
D\to\partial\mathbb D$ be a function of bounded variation and
$\varphi:\partial D\to\mathbb R$ be a measurable function with
respect to logarithmic capacity. Then there is an analytic function
$f:D\to\Bbb C$ such that
\begin{equation}\label{26} \lim\limits_{z\to\zeta}\ \mathrm {Re}\
\{\overline{\lambda(\zeta)}\cdot f(z)\}\ =\ \varphi(\zeta)
\quad\quad\quad\mbox{for}\ \ \mbox{a.e.}\ \ \ \zeta\in\partial D
\end{equation}
with respect to logarithmic capacity in the sense of the unique
principal asymptotic value.}
\end{theorem}

\medskip

In particular, choosing $\lambda\equiv1$ in~\eqref{26}, we obtain
the following consequence.

\medskip

\begin{proposition}\label{P3} {\it
Let $D$ be a Jordan domain and let $\varphi:\partial D\to\mathbb R$
be a measurable function with respect to logarithmic capacity. Then
there is an analytic function $f:D\to\Bbb C$ such that
\begin{equation}\label{27}
\lim\limits_{z\to\zeta}\ \mathrm {Re}\ f(z)\ =\ \varphi(\zeta)
\quad\quad\quad\mbox{for}\ \ \mbox{a.e.}\ \ \ \zeta\in\partial D
\end{equation}
with respect to logarithmic capacity in the sense of the unique
principal asymptotic value.}
\end{proposition}

\medskip

\begin{corollary}\label{C2} {\it
Under the conditions of proposition \ref{P3}, there is a harmonic
function $u$ in $D$ such that in the same sense
\begin{equation}\label{19} \lim\limits_{z\to\zeta}\ \mathrm
u(z)\ =\ \varphi(\zeta) \quad\quad\quad\mbox{for}\ \ \mbox{a.e.}\ \
\ \zeta\in\partial D\ .
\end{equation}}
\end{corollary}

\medskip

\begin{remark}\label{R5}
It is easy to see that here in comparison with the paper \cite{Ry},
we strengthen the conditions as well as the conclusions of these
theorems, see Remark \ref{R1}.
\end{remark}

\bigskip

\section{The Riemann-Hilbert problem in quasidisks}

\medskip

\begin{theorem}\label{T6} {\it
Let $D$ be a Jordan domain in $\Bbb C$ bounded by a quasiconformal
curve, $\mu:D\to\mathbb{C}$ be a measurable (by Lebesgue) function
with $||\mu||_{\infty}<1$, $\lambda:\partial D\to\mathbb{C},\:
|\lambda(\zeta)|\equiv1$, be a function of bounded variation and let
$\varphi:\partial D\to\mathbb{R}$ be a measurable function with
respect to logarithmic capacity. Then the Beltrami equation
\eqref{1} has a regular solution of the Riemann-Hilbert problem
\eqref{2}. If in addition $\partial D$ is rectifiable, then the
limit in \eqref{2} holds a.e. with respect to the natural parameter
along any nontangential path.}
\end{theorem}

\bigskip

In particular, the latter conclusion of Theorem \ref{T6} holds in
the case of smooth boundaries.

\bigskip

\begin{proof}
Without loss of generality we may assume that $0\in D$ and
$1\in\partial D$. Extending $\mu$ by zero everywhere outside of $D$,
we obtain the existence of a quasiconformal mapping
$f:\overline{\mathbb{C}}\to\overline{\mathbb{C}}$ with the
normalization $f(0)=0,\: f(1)=1$ и $f(\infty)=\infty$ satisfying the
Beltrami equation \eqref{1} with the given $\mu$, see e.g. Theorem
V.B.3 in \cite{Alf}. By the theorems of Riemann and Caratheodory,
the Jordan domain $f(D)$ can be mapped by a conformal mapping $g$
with the normalization $g(0)=0$ and $g(1)=1$ onto the unit disk
$\Bbb D$. It is clear that $h:=g\circ f$ is a quasiconformal
homeomorphism with normalization $h(0)=0$ и $h(1)=1$ satisfying the
same Beltrami equation.

By the reflection principle for quasiconformal mappings, using the
conformal reflection (inversion) with respect to the unit circle in
the image and qua\-si\-con\-for\-mal reflection with respect to
$\partial D$ in the preimage, we can extend $h$ to a
qua\-si\-con\-for\-mal mapping
$H:\overline{\mathbb{C}}\to\overline{\mathbb{C}}$ with the
normalization $H(0)=0,\: H(1)=1$ и $H(\infty)=\infty$, see e.g.
I.8.4, II.8.2 and II.8.3 in \cite{LV}. Note that
$\Lambda=\lambda\circ H^{-1}$ is a function of bounded variation,
$V_{\Lambda}(\partial{\mathbb{D}})=V_{\lambda}(\partial{{D}})$.

The mappings $H$ and $H^{-1}$ transform sets of logarithmic capacity
zero on $\partial D$ into sets of logarithmic capacity zero on
$\partial\Bbb D$ and vice versa because quasiconformal mappings are
continuous by H\"older on $\partial D$ and $\partial\Bbb D$
correspondingly, see e.g. Theorem II.4.3 in \cite{LV}.

Further, the function $\Phi=\varphi\circ H^{-1}$ is measurable with
respect to logarithmic capacity. Indeed, under this mapping
measurable sets with respect to lo\-ga\-rith\-mic capacity are
transformed into measurable sets with respect to logarithmic
capacity because such a set can be represented as the union of a
sigma-com\-pac\-tum and a set of logarithmic capacity zero and
compacta under continuous mappings are transformed into compacta and
compacta are measurable sets with respect to logarithmic capacity.

Thus, the original Riemann-Hilbert problem for the Beltrami equation
\eqref{1} is reduced to the Riemann-Hilbert problem for analytic
functions $F$ in the unit circle:
\begin{equation}\label{29}
\lim_{z\to\zeta}\ \overline{\Lambda(\zeta)}\cdot F(z)\ =\
\Phi(\zeta)
\end{equation}
and by Theorem \ref{T4} there is an analytic function $F:\mathbb
D\to\mathbb C$ for which this boundary condition holds for a.e.
$\zeta\in\partial\mathbb D$ with respect to logarithmic capacity
along any nontangential path.

\medskip

So, the desired solution of the original Riemann-Hilbert problem
\eqref{2} for the Beltrami equation \eqref{1} exists and can be
represented as $f=F\circ H$.

\medskip

Finally, since the distortion of angles under the quasiconformal
mapping is bounded, see e.g. \cite{A}--\cite{Ta}, then in the case
of a rectifiable boundary of $D$ condition \eqref{2} can be
understood along any nontangential path a.e. with respect to the
natural parameter.
\end{proof}$\Box$

\bigskip

\section{On the dimension of spaces of solutions}

\medskip

By the known Lindel\"of maximum principle, see e.g. Lemma 1.1 in
\cite{GM}, it follows the uniqueness theorem for the Dirichlet
problem in the class of {\bf bounded} harmonic functions on the unit
disk $\mathbb D = \{ z\in\mathbb{C}: |z|<1\}$. In general there is
no uniqueness theorem in the Dirichlet problem for the Laplace
equation. Furthermore, it was proved in \cite{Ry} that the space of
all harmonic functions in $\mathbb D$ with nontangential limit $0$
at a.e. point of $\partial\mathbb D$ has the infinite dimension.

\medskip

Let us show that in view of Lemma 3.1 one can similarly prove the
more refined results on harmonic functions with respect to
logarithmic capacity instead of the measure of the length on
$\partial\mathbb D$.

\bigskip

\begin{theorem}\label{T7}
{\it The space of all harmonic functions $u:\mathbb D\to\mathbb R$
such that $\lim\limits_{z\to\zeta}u(z) = 0$ along any nontangential
path for a.e. $\zeta\in\partial\mathbb D$ with respect to
logarithmic capacity has the infinite dimension.}
\end{theorem}

\medskip

\begin{proof} Indeed, let $\Phi:[0,2\pi]\to\mathbb R$ be integrable, differentiable
and $\Phi^{\prime}(t)=0$  a.e. on $\partial\mathbb D$ with respect
to logarithmic capacity. Then the function
$$
U(z)\ \colon =\ \frac{1}{2\pi}\
\int\limits_{0}\limits^{2\pi}\frac{1-r^2}{1-2r\cos(\vartheta-t)+r^2}\
\Phi(t)\ dt\ ,\ \ \ z=re^{i\vartheta}, \ r<1\ ,
$$
is harmonic on $\mathbb D$ with $U(z)\to\Phi(\Theta)$ as $z\to
e^{i\Theta}$, see e.g. Theorem 1.3 in \cite{GM} or Theorem IX.1.1 in
\cite{Go}, and $\frac{\partial}{\partial\vartheta}\ U(z)\to
\Phi^{\prime}(\Theta)$ as $z\to e^{i\Theta}$ along any nontangential
path whenever $\Phi^{\prime}(\Theta)$ exists, see e.g. 3.441 in
\cite{Z}, p. 53, or Theorem IX.1.2 in \cite{Go}. Thus, the harmonic
function $u(z)\ =\ \frac{\partial}{\partial\vartheta}\ U(z)$ has
nontangential limit 0 at a.e. point of $\partial\mathbb D$ with
respect to logarithmic capacity.

Let us give a subspace of such functions $u$ with an infinite basis.
Namely, let $\varphi:[0,1]\to[0,1]$ be a function of  the Cantor
type, see Lemma \ref{L2}, and let $\varphi_n:[0,2\pi]\to[0,1]$ be
equal to $\varphi((t-a_{n-1})/(a_n-a_{n-1}))$ on $[a_{n-1},a_n)$
where $a_0=0$ and $a_n=2\pi(2^{-1}+\ldots +2^{-n})$, $n=1,2,\ldots $
and $0$ outside of $[a_{n-1},a_n)$. Denote by $U_n$ and $u_n$ the
harmonic functions corresponding to $\varphi_n$ as in the first
item.

By the construction the supports of the functions $\varphi_n$ are
mutually disjoint and, thus, the series
$\sum\limits_{n=1}\limits^{\infty}\gamma_n\varphi_n$ is well defined
for every sequence $\gamma_n\in\mathbb R$, $n=1,2,\ldots $. If in
addition we restrict ourselves to the sequences
$\gamma=\{\gamma_n\}$ in the space $l$ with the norm
$\|\gamma\|=\sum\limits_{n=1}\limits^{\infty}|\gamma_n|<\infty$,
then the series is a suitable function $\Phi$ for the first item.

Denote by $U$ and $u$ the harmonic functions corresponding to the
function $\Phi$ as in the first item and by ${\mathcal H}_0$ the
class of all such $u$. Note that $u_n$, $n=1,2,\ldots $,  form a
basis in the space ${\mathcal H}_0$ with the locally uniform
convergence in $\mathbb D$ which is metrizable.

\medskip

Firstly, $\sum\limits_{n=1}\limits^{\infty}\gamma_nu_n\ne 0$ if
$\gamma\ne 0$. Really, let us assume that $\gamma_n\ne 0$ for some
$n=1,2,\ldots $. Then $u\ne 0$ because the limits
$\lim\limits_{z\to\zeta}U(z)$ exist for all $\zeta=e^{i\vartheta}$
with $\vartheta\in(a_{n-1},a_n)$ and can be arbitrarily close to $0$
as well as to $\gamma_n$.

\bigskip

Secondly, $u^*_m=\sum\limits_{n=1}\limits^{m}\gamma_nu_n\to u$
locally uniformly in $\mathbb D$ as $m\to\infty$. Indeed, elementary
calculations give the following estimate of the remainder term
$$
|u(z)-u^*_m(z)|\ \leq\ \frac{2r(1+r)}{(1-r)^3}\ \cdot
\sum\limits_{n=m+1}\limits^{\infty}|\gamma_n|\ \to\ 0\ \ \ \ \ \ \ \
\ {\mbox{az}}\ \ \ \ \ m\to\infty
$$
in every disk $\mathbb D(r)=\{ z\in\mathbb C: |z|\leq r\}$, $r<1$.
\end{proof}

\medskip

\begin{remark}\label{INT} Note that harmonic functions $u$ found by us in Theorem \ref{T7}
themselves cannot be represented in the form of the Poisson integral
with any function $\Phi\in\, L_p([0,2\pi])$, $p>1$, because this
integral would have nontangential limits $\Phi$ a.e., see e.g.
Corollary IX.9.1 in \cite{Go}. Thus, $u$ do not belong to the
classes $h_p$ for any $p>1$, see e.g. Theorem IX.2.3 in \cite{Go}.
\end{remark}

\bigskip

\begin{corollary}\label{C3} {\it Given a measurable function $
\varphi :\partial\mathbb D\to \mathbb R$, the space of all harmonic
functions $u:\mathbb D\to\mathbb R$ with the limits
$\lim\limits_{z\to\zeta}u(z) = \varphi(\zeta)$ for a.e.
$\zeta\in\partial\mathbb D$ with respect to logarithmic capcity
along nontangential paths has the infinite dimension.}
\end{corollary}

\medskip

Indeed, we have at least one such harmonic function $u$ by Theorem
\ref{T2} and, combining this fact with Theorem \ref{T7}, we obtain
the conclusion of Corollary \ref{C3}.

\bigskip

The statements on the infinite dimension of the space of solutions
can be extended to the Riemann-Hilbert problem because we have
reduced this problem in Theorem \ref{T4}  to the corresponding two
Dirichlet problems.

\medskip

\begin{theorem}\label{T8} {\it
Let $\lambda:\partial\mathbb D\to\partial\mathbb D$ be of bounded
variation and $\varphi:\partial\mathbb D\to\mathbb R$ be measurable
with respect to logarithmic capacity. Then the space of all analytic
functions $f:\mathbb D\to\mathbb C$ such that along any
nontangential path
\begin{equation}\label{25}
\lim\limits_{z\to\zeta}\ \mathrm {Re}\
\{\overline{\lambda(\zeta)}\cdot f(z)\}\ =\ \varphi(\zeta)
\quad\quad\quad\mbox{for}\ \ \mbox{a.e.}\ \ \
\zeta\in\partial\mathbb D
\end{equation}
with respect to logarithmic capacity has the infinite dimension.}
\end{theorem}

\medskip

\begin{proof} Let $u:\mathbb D\to\mathbb R$ be a harmonic function with nontangential limit $0$ at
a.e. point of $\partial\mathbb D$ with respect to logarithmic
capacity from Theorem \ref{T7}. Then there is the unique harmonic
function $v:\mathbb D\to\mathbb R$ with $v(0)=0$ such that
${\mathcal C}=u+iv$ is an analytic function. Thus, setting in the
proof of Theorem \ref{T4} $g={\mathcal A}({\mathcal B}+{\mathcal
C})$ instead of $f={\mathcal A}\cdot{\mathcal B}$, we obtain by
Theorem \ref{T7} the space of solutions of the Riemann-Hilbert
problem (\ref{25}) for analytic functions of the infinite dimension.
\end{proof}

\bigskip

\begin{remark}\label{R1DIM} The dimension of the spaces of solutions of
the Riemann-Hilbert problem for the Beltrami equation in Theorem
\ref{T6} is also infinite because this case is reduced  to the case
of Theorem \ref{T8} as in the proof of Theorem \ref{T6}.
\end{remark}

\bigskip

\section{Extension of results to countably bounded variation}

\medskip

We call $\lambda:\partial\mathbb D\to\mathbb C$ a function of {\bf
countably bounded variation}, write
$\lambda\in\mathcal{CBV}(\partial\mathbb D)$, if there is a
countable collection of mutually disjoint arcs $\gamma_n$,
$n=1,2,\ldots$ on each of which the restriction of $\lambda$ is of
bounded variation $V_n$, $\sum\limits_{n=1}\limits^{\infty} V_n\cdot
|\gamma_n|<\infty$, and the set $\partial\mathbb D\setminus
\bigcup\limits_1\limits^{\infty}\gamma_n$ is countable. The
definition is also extended in a natural way to an arbitrary Jordan
curve $\Gamma$ in $\mathbb C$. All the above results on the
Riemann-Hilbert problem have been extended to the case of
$\lambda\in\mathcal{CBV}$ in the paper \cite{Ye}. The latter was
based on the following analogs of Proposition \ref{P2} and Lemma
\ref{T3}.

\bigskip

\begin{proposition}\label{P22} {\it
For every function $\lambda:\partial\mathbb{D}\to\partial\mathbb{D}$
of the class $\mathcal{CBV}(\partial\mathbb{D})$ there is a function
$\alpha_{\lambda}:\partial\mathbb{D}\to\mathbb{R}$ of the class
$L_1(\partial\mathbb D)\cap\mathcal{CBV}(\partial\mathbb{D})$ such
that $\lambda(\zeta)=\exp\{{i\alpha_{\lambda}}(\zeta)\}$,
$\zeta\in\partial\mathbb{D}$.}
\end{proposition}

\bigskip

\begin{proof}
Denote by $\lambda_n$ the complex valued function on
$\partial\mathbb{D}$ that is equal to $\lambda$ on $\gamma_n$ and to
$1$ outside of $\gamma_n$. Let $\alpha_n$ correspond to $\lambda_n$
by Proposition \ref{P2}. Then its variation $V_n^*\le V_n\cdot
3\pi/2$. With no loss of generality we may assume that
$\alpha_n\equiv 0$ outside of $\gamma_n$. Set $\alpha
=\sum\limits_{n=1}\limits^{\infty}\alpha_n$. Then
$\alpha\in\mathcal{CBV}(\partial\mathbb{D})$ and
$\lambda(\zeta)=\exp\{{i\alpha}(\zeta)\}$,
$\zeta\in\partial\mathbb{D}$. Applying the corresponding shifts
(divisible $2\pi$) we may change $\alpha_n$ on $\gamma_n$ through
$\alpha_n^*$ with $|\alpha_n^*|\le\pi$ at the middle point of
$\gamma_n$. Then it is clear that the new function
$\alpha^*\in\mathcal{CBV}(\partial\mathbb{D})$ and
$\lambda(\zeta)=\exp\{{i\alpha^*}(\zeta)\}$,
$\zeta\in\partial\mathbb{D}$, and, moreover, $|\alpha^*|\le\pi +
V_n\cdot 3\pi/2$ on every $\gamma_n$ and hence $\| \alpha^*\|_1\le
2\pi^2+\frac{3\pi}{2}\sum\limits_{n=1}\limits^{\infty}V_n<\infty$,
i.e. $\alpha^*\in L_1(\partial\mathbb D)$. \end{proof}$\Box$

\bigskip

We prove the following statement similarly to Lemma \ref{T3} ,
however, without the reference to the paper \cite{T} because in the
case the function $\alpha$, generally speaking, will be not of
bounded variation but we are able instead it to apply Lemma \ref{T3}
itself.

\bigskip

\begin{lemma}\label{T3C} {\it
Let $\mathbb D$ be the unit disk in $\mathbb C$,
$\alpha:\partial\mathbb D\to\mathbb R$ be a bounded  function of the
class $\mathcal{CBV}(\partial\mathbb{D})$, $u:\mathbb D\to\mathbb R$
be a bounded harmonic function such that
\begin{equation}\label{23C}
\lim\limits_{z\to\zeta}\ u(z)\ =\ \alpha(\zeta)
\end{equation}
at every point of continuity of $\alpha$ and let $v$ be its
conjugate harmonic function. Then for a.e. $\zeta\in\partial\mathbb
D$ with respect to logarithmic capacity
\begin{equation}\label{24C}
\lim\limits_{z\to\zeta}\ v(z)\ =\ \beta(\zeta)
\end{equation}
along any nontangential path in $\mathbb D$ terminating at $\zeta$
where $\beta:\partial\mathbb D\to\mathbb R$ is a function that is
measurable with respect to logarithmic capacity.}
\end{lemma}

\medskip

\begin{proof}
Indeed, the function $\alpha\in\mathcal{CBV}(\partial\mathbb{D})$
has at most a countable set $S$ of points of discontinuity. Hence by
the generalized maximum principle such a function $u$ is unique.
Moreover, $\alpha\in L_1(\partial\mathbb D)$ and, consequently, $u$
can be represented as the Poisson integral of the function $\alpha$
\begin{equation}\label{PoissonVC}
u(re^{i\vartheta})\ =\ \frac{1}{2\pi}\
\int\limits_{-\pi}\limits^{\pi}\frac{1-r^2}{1-2r\cos(\vartheta
-t)+r^2}\ \alpha(e^{it})\ dt\ .
\end{equation}
The Schwartz integral
\begin{equation}\label{PoissonWC}
f(z)\ :=\ \frac{1}{2\pi i}\ \int\limits_{\partial\mathbb D}
\alpha(\zeta)\ \frac{\zeta +z}{\zeta -z}\ \frac{d\zeta}{\zeta}
\end{equation}
gives the analytic function $f=u+iv$ in $\mathbb D$, where
\begin{equation}\label{PoissonVC}
v(re^{i\vartheta})\ =\ \frac{1}{2\pi}\
\int\limits_{-\pi}\limits^{\pi}\frac{2r\sin(\vartheta-t)}{1-2r\cos(\vartheta-t)+r^2}\
\alpha(e^{it})\ dt\ .
\end{equation}

Let us apply the linearity of the integral operator
(\ref{PoissonVC}). Namely, denote by $\chi$ the characteristic
function of an arc $\gamma_*$ of $\partial\mathbb D$ where $\alpha$
is of bounded variation from the definition of  $\mathcal{CBV}$.
Setting $\alpha_*=\alpha\cdot\chi$ and $\alpha_0=\alpha -\alpha_*$,
we have that $\alpha =\alpha_*+\alpha_0$. Then $v=v_*+v_0$ where
$v_*$ and $v_0$ correspond to $\alpha_*$ and $\alpha_0$ by formula
(\ref{PoissonVC}). By Lemma \ref{T3} for a.e.
$\zeta\in\partial\mathbb D$ with respect to logarithmic capacity
\begin{equation}\label{24CC}
\lim\limits_{z\to\zeta}\ v_*(z)\ =\ \beta_*(\zeta)
\end{equation}
along any nontangential path in $\mathbb D$ terminating at $\zeta$
where $\beta_*:\partial\mathbb D\to\mathbb R$ is a function that is
measurable with respect to logarithmic capacity. Moreover, it is
evident from formula (\ref{PoissonVC}) that $v_0(z)\to
\beta_0(\zeta)$ as $z\to\zeta$ for all $\zeta\in\gamma_*$ where
$\beta_0:\gamma_*\to\mathbb R$ is even continuous on $\gamma_*$.
Thus, setting $\beta=\beta_*+\beta_0$ on $\gamma_*$, we obtain the
conclusion of Lemma \ref{T3C} by countability of the collection of
such arcs $\gamma_*$ and by countability of the completion of this
collection on $\partial\mathbb D$.
\end{proof}$\Box$

\medskip
\noindent
{\bf Artyem Efimushkin and Vladimir Ryazanov,}\\
Institute of Applied Mathematics and Mechanics,\\
National Academy of Sciences of Ukraine,\\
74 Roze Luxemburg Str., Donetsk, 83114,\\
art89@bk.ru, vl.ryazanov1@gmail.com

\end{document}